\documentclass[12pt]{amsart}
\usepackage{amssymb}
\usepackage{hyperref}

\theoremstyle{plain}
\newtheorem{theorem}{Theorem}

\newtheorem{corollary}{Corollary}
\newtheorem{proposition}{Proposition}
\newtheorem{definition}{Definition}

\newtheorem*{remark}{Remark}

\title[$b$-additive and $b$-multiplicative Ramanujan-Hardy numbers]{Infinite sets of $b$-additive and $b$-multiplicative Ramanujan-Hardy numbers}

\author{Viorel Ni\c{t}ic\u{a}}
\address{Department of Mathematics\\ West Chester
University of Pennsylvania\\ West Chester, PA 19383, USA\\
vnitica@wcupa.edu}

\begin{document}

\begin{abstract} Let $b$ be a numeration base. A $b$-additive Ramanujan-Hardy number $N$ is an integer for which there exists at least an integer $M$, called additive multiplier, such that the product of $M$ and the sum of base $b$ digits of $N$,  added to the reversal of the product, gives $N$. We show that for any $b$ there exists an infinity of $b$-additive Ramanujan-Hardy numbers and an infinity of additive multipliers. A $b$-multiplicative Ramanujan-Hardy number $N$ is an integer for which there exists at least an integer $M$, called multiplicative multiplier, such that the product of $M$ and the sum of base $b$ digits of $N$,  multiplied by the reversal of the product, gives $N$. We show that for an even $b$, $b\equiv 1 \pmod {3}$, and for $b=2$, there exists an infinity of $b$-multiplicative Ramanujan-Hardy numbers and an infinity of multiplicative multipliers. 

These results completely answer two questions and partially answer two other questions among those asked in V. Ni\c tic\u a, \emph{About some relatives of the taxicab number}, arXiv:1805.10739v3.
\end{abstract}

\maketitle

\section{Introduction}\label{sec:1}

Let $b\ge 2$ be a numeration base. In Ni\c tic\u a \cite{N}, motivated by some properties of the taxicab number, 1729, we introduce the classes of $b$-additive Ramanujan-Hardy (or $b$-ARH) numbers and $b$-multiplicative Ramanujan-Hardy (or $b$-MRH) numbers. The first class consists of numbers for which there exists at least an integer $M$, called \emph{additive multiplier}, such that the product of $M$ and the sum of base $b$ digits of $N$,  added to the reversal of the product, gives $N$. The second class consists of numbers for which there exists at least an integer $M$, called \emph{multiplicative multiplier}, such that the product of $M$ and the sum of base $b$ digits of $N$,  multiplied by the reversal of the product, gives $N$. 

It is asked \cite[Question 6]{N} if the set of $b$-ARH numbers is infinite and it is asked \cite[Question 8]{N} if the set of additive multipliers is infinite. It is shown \cite[Theorems 12 and 15]{N} that the answer is positive if $b$ is even. The case $b$ odd is left open. It is asked \cite[Question 7]{N} if the set of $b$-MRH numbers is infinite for all numeration  bases and it is asked \cite[Question 9]{N} if the set of multiplicative multipliers is infinite. It is shown \cite[Theorem 30]{N} that the answer is positive if $b$ is odd. The case $b$ even is left open. 

We recall that \emph{Niven (or Harshad) numbers} are numbers divisible by the sum of their decimal digits. Niven numbers have been extensively studied. See for instance Cai \cite{C}, Cooper and Kennedy \cite{CK}, De Koninck and Doyon \cite{KD}, and Grundman \cite{G}. Of interest are also $b$-Niven numbers, which are numbers divisible by the sum of their base $b$ digits. See, for example, Fredricksen, Iona\c scu, Luca, and St\u anic\u a \cite{FILS}. A $b$-MRH-number is a $b$-Niven number. Nigh degree $b$-Niven numbers are introduced in \cite{N2}.

The goal of this paper is to show that, for any numeration base, there exist an infinity of $b$-ARH numbers and an infinity of distinct additive multipliers. We also show that for even $b$, $b\equiv 1 \pmod {3}$, and for $b=2$,  there exist an infinity of $b$-MRH numbers, and an infinity of distinct multiplicative multipliers. The results here overlap with some in \cite{N}, but with different sets of examples. They also completely answer the first two questions from \cite{N} revisited above, and partially answer the other two. We observe that a trivial example of an infinity of $b$-MRH numbers is given by $\{[1(0)^{\land k}]_b\vert k\in \mathbb{N}\}$. The examples we show here have at least two digits different from zero. Finding an infinity of $b$-MARH numbers with all digits different from zero remains an open question.

Our results about $b$-ARH numbers also give solutions to the dyophantine equation $N\cdot M=reversal(N\cdot M)$. Motivated by this link, we show that the dyophantine equation has a solution for all fixed integers $N$ not divisible by the bse $b$ and for any nymeration base. Our final result shows that for any string of digits $I$ there exists an infinity of $b$-Niven number that contains $I$ in their base $b$-representation. We do not know how to prove a similar result for the classes of $b$-ARH and $b$-MRH numbers. 

\section{Statements of the main results}\label{sec:1-main}

Let $s_b(N)$ denote the sum of base $b$ digits of integer $N$. If $x$ is a string of digits and let $(x)^{\land k}$ denoye the base $10$ integer obtained by repeating $x$ $k$-times. Let $[x]_b$ denote the value of the string $x$ in base $b$. The \emph{reversal} of an integer $N$ is the number obtained from $N$ writing its digits in reverse order. If $N$ is an integer written in base $b$, let $N^R$ denote the reversal of $N$. While the operations of addition and multiplication of integers are independent of the base, the operation 
of taking the reversal is not. In the definition of a $b$-ARH-number/$b$-MRH number $N$ we take the reversal of the base $b$ representation of $s_b(N)M$.

\begin{theorem}\label{thm:1} Let $\alpha\ge 1$ integer, $b\ge \alpha+1$ integer, and $k=(1+\alpha)^\ell, \ell\ge 0$.
Assume $b\equiv 2+\alpha \pmod{2+2\alpha}$. Define
\begin{equation*}
N_k=[(1\alpha)^{\land k}]_b.
\end{equation*}

Then there exists $M\ge 0$ integer such that
\begin{equation*}
s_b(N_k)\cdot M=(s_b(N_k)\cdot M)^R=\frac{N_k}{2}.
\end{equation*}

In particular, the numbers $N_k, k\ge 1,$ are $b$-ARH numbers and $b$-Niven numbers.
\end{theorem}

The proof of Theorem \ref{thm:1} is done in Section \ref{sec:2}.

\begin{remark} The particular case $b=10, \alpha=2,$ of Theorem \ref{thm:1}, which gives $N_k=(12)^{3^\ell}$, is covered by \cite[Example 10]{N}.
Theorem \ref{thm:1} does not give any information if $b=2$. 
\end{remark}

The following proposition gives positive answers to \cite[Questions 5 and 6]{N}.

\begin{proposition}\label{prop:1} For any numeration base, there exists an infinite set of $b$-ARH numbers and an infinite set of additive multipliers. The $b$-ARH numbers in the infinite set also are $b$-Niven numbers.
\end{proposition}

The proof of Proposition \ref{prop:1} is done in Section \ref{sec:3-pre}.

\begin{remark} We observe that \cite[Theorems 12 and 15]{N} show, for all even bases, an infinity of $b$-ARH numbers that are not $b$-Niven numbers. The case of odd base is open. The question of finding an infinity of $b$-Niven numbers that are not $b$-ARH numbers is also open. It is shown in \cite[Theorem 28]{N} that for any base there exists an infinity of numbers that are not $b$-ARH numbers.
\end{remark}

The result in Theorem \ref{thm:1} gives many base $10$ solutions for the equation:
\begin{equation}\label{eq:star10}
N\cdot M=(N\cdot M)^R.
\end{equation}

One can try to solve, for any numeration base $b$, the equation:
\begin{equation}\label{eq:star}
N\cdot M=(N\cdot M)^R,
\end{equation}
where $(N\cdot M)^R$ is the reversal of $N\cdot M$ written in base $b$.

Observe that if $N$ is divisible by $b$, then $(N\cdot_b M)^R$ has less digits then $N\cdot_b M$, therefore $N$ is not a solution of \eqref{eq:star}. Note also that if $N=N^R$ and $N$ has $k$ digits then  \eqref{eq:star} always has an infinite set of solutions with 
\begin{equation*}
M=[(1(0)^{\land\ell})^{\land p}1]_b, \ell \ge k-1, p\ge 0.
\end{equation*}
Consequently, if $(N_0,M_0)$ is a solution of \eqref{eq:star}, then \eqref{eq:star} also has infinite sets of solutions of types $(N_0,M)$ and $(N,M_0)$.

\begin{theorem}\label{thm:new} Let $b\ge 2$ and $N\ge 1$ integer such that $b\not \vert N$. Then $N$ is a solution of \eqref{eq:star}. 
\end{theorem}

 The proof of Theorem \ref{thm:new} is done in Section \ref{sec:3}. For base 10, a proof belonging to David Radcliffe can be found at   \cite{Rat}. We learned about the reference \cite{Rat} from J. Shallit. We generalize the proof for an arbitrary numeration base. After this paper was written, we learned from J. Shallit \cite{S} that he also has a proof of Theorem \ref{thm:new}. 

A \emph{$b$-numeric palindrome} is a base $b$ integer $N$ such that $N=N^R$.

\begin{corollary} The prime factors of $b$-numeric palindromes exhaust the set of prime numbers. 
\end{corollary}

\begin{definition} The \emph{multiplicity} of a multiplier $M$ is the number of $N$ solutions of  \eqref{eq:star}.
\end{definition}

It was observed above that for any solution $(N,M)$ of \eqref{eq:star}, $M$ has infinite multiplicity. The following theorem shows infinite sets of solutions of \eqref{eq:star} that cannot be derived from the previous considerations.

\begin{theorem}\label{thm:33} Let $b\ge 2$ a numeration base. Then:
\begin{equation*}
[1(b-1)]_b\cdot [(b-1)^{\land k}]_b=[1(b-2)(b-1)^{\land k-2}(b-2)1]_b
\end{equation*}
for all $k\ge 0$. 
\end{theorem}

The proof of Theorem \ref{thm:33} is done in Section \ref{sec:44}.

In \cite{N4} we show more infinite sets of solutions of \eqref{eq:star}. We also show infinite sets of numbers that satisfy equation \eqref{eq:star} up to a prescribed number of misplaced digits for which we know their position. 

Our next results shows, for $b$ even, more examples of infinite sets of $b$-ARH that are not $b$-Niven numbers.

\begin{theorem}\label{thm:final} Let $b\ge 2$ even. Let $a\in\{1,2,\dots,b-1\}$ and $k\ge 0$ integer. 

a) Let
\begin{equation*}
N_k=[a(0)^{\land k}a]_b.
\end{equation*}
Then $N_k$ is a $b$-ARH number, but not a $b$-Niven number.

b) Let
\begin{equation*}
N_k=[\left ( 1(0)^{\land k} \right )^{\land b}0 \left ( (0)^k 1\right )^{\land b}]_b.
\end{equation*}
Then $N_k$ is a $b$-ARH number, but not a $b$-Niven number

c) Let
\begin{equation*}
N_k=[\left ( (0)^{\land k}1 \right )^{\land b}0 \left (1 (0)^k \right )^{\land b}]_b.
\end{equation*}
Then $N_k$ is a $b$-ARH number and  a $b$-Niven number.
\end{theorem}

The proof of Theorem \ref{thm:final} is done in Section \ref{sec:final}.

\begin{theorem}\label{thm:final2} a) Let $b\ge 4$ even and $b\equiv 1\pmod{3}$. Let $k\ge 1$ integer such that $k\equiv 1\pmod{3}$. Define
\begin{equation*}
\begin{gathered}
\alpha_k=[1(0)^{\land k}(b-2)]_b.
\end{gathered}
\end{equation*}

Then $N_k=\alpha_k\cdot(\alpha_k)^R$ is a $b$-MRH number.

b) Let $b=2$ and $k\ge 1$ even integer. Define
\begin{equation*}
\alpha_k=[1(0)^k1]_2.
\end{equation*}

Then $N_k=\alpha_k\cdot(\alpha_k)^R$ is a $b$-MRH number.
\end{theorem}

The proof of Theorem \ref{thm:final2} is done in Section \ref{sec:final2}.

The following proposition gives partial answers to \cite[Questions 7 and 8]{N}.

\begin{proposition}\label{prop:1-mult} For any even numeration base $b$, $b\equiv 1\pmod{3}$ and for $b=2$ there exist an infinite set of $b$-MRH numbers and an infinite set of multipliers.
\end{proposition}

The proof of Proposition \ref{prop:1-mult} is done in Section \ref{sec:3-pre-mult}.

Our next result lists several infinite sequences of $10$-MRH-numbers.

\begin{proposition} Assume $k\ge 1$ integer and define $N_k=\alpha_k\cdot (\alpha_k)^R$, where $\alpha_k$ is one of the following numbers:
\begin{itemize}
\item $[1(0)^{\land k} 8]_{10}, k\equiv 1\pmod{3}$,
\item $[7(0)^{\land k} 2]_{10},$
\item $[5(0)^{\land k} 4]_{10},$
\item $[4(0)^{\land k} 5]_{10}$
\end{itemize}
Then $N_k$ is an 10-MRH number.
\end{proposition}

The first item follows an a corollary of Theorem \ref{thm:final2}. The other items can be proved using the same approach as in the proof of Theorem \ref{thm:final2}.

\begin{theorem}\label{thm:last-thm} For any base $b$ and for any string of base $b$ digits $I$ there exists an infinity of $b$-Niven number that contains the string $I$ in their base $b$-representation.
\end{theorem}

The proof of Theorem \ref{thm:last-thm} is done in section \ref{seq:last-thm}.

\section{Proof of Theorem \ref{thm:1}}\label{sec:2}

\begin{proof} The base $b$ representation for $N_k/2$ is $N_k/2=\left [ \left ( 0\frac{b+\alpha}{2}\right )^{\land k}\right ]_b$. One has that:
\begin{equation}\label{eq:55}
s_b(N_k)=k\cdot(1+\alpha)=(1+\alpha)^{\ell+1}.
\end{equation}

The value of $N_k/2$ in base 10 is obtained summing a geometric series.
\begin{equation}\label{eq:66}
\begin{gathered}
\frac{N_k}{2}=\frac{b+\alpha}{2}\cdot b^{2k-2}+\frac{b+\alpha}{2}\cdot b^{2k-4}+\cdots\\
 + \frac{b+\alpha}{2}\cdot b^{2}+\frac{b+\alpha}{2}=\frac{b+\alpha}{2}\cdot \frac{b^{2k}-1}{b^2-1}\\
=\frac{b+\alpha}{2}\cdot \frac{b^{2(1+\alpha)^{\ell}}-1}{b^2-1}.
\end{gathered}
\end{equation}

Note that $N_k/2=(N_k/2)^R$. We finish the proof of the theorem if we show that:
\begin{equation}\label{eq:1}
(1+\alpha)^{\ell+1} \Big \vert \frac{b+\alpha}{2}\cdot \frac{b^{2(1+\alpha)^{\ell}}-1}{b^2-1}.
\end{equation} 

We prove \eqref{eq:1} by induction on $\ell$. For $\ell=0$ equation \eqref{eq:1} becomes
\begin{equation*}
1+\alpha \Big \vert \frac{b+\alpha}{2},
\end{equation*}
which is true because $b\equiv 2+\alpha \pmod{2+2\alpha}$.

Now we assume that \eqref{eq:1} is true for $\ell$ and show that it is true for $\ell+1$.  

\begin{equation}\label{eq:account}
\begin{gathered}
\frac{b+\alpha}{2}\cdot \frac{b^{2(1+\alpha)^{\ell+1}}-1}{b^2-1}=\frac{b+\alpha}{2}\cdot \frac{\left ( b^{2(1+\alpha)^{\ell}}\right )^{1+\alpha}-1}{b^2-1}\\
=\frac{b+\alpha}{2}\cdot \frac{b^{2(1+\alpha)^{\ell}}-1}{b^2-1}\left ( B^{\alpha}+B^{\alpha-1}+\cdots +B^2+B+1\right ),
\end{gathered}
\end{equation}
where
\begin{equation}\label{eq:B}
B=b^{2(1+\alpha)^{\ell}}.
\end{equation}

The congruence $b\equiv 2+\alpha \pmod{2+2\alpha}$ implies that
\begin{equation*}
b^2\equiv (2+\alpha)^2\equiv \alpha^2+4\alpha+4\equiv \alpha^2\equiv 1 \pmod{1+\alpha},
\end{equation*}
which implies that
\begin{equation}\label{eq:bsquare}
b^m\equiv 1 \pmod{1+\alpha}, \text{ $m$ even}.
\end{equation}

From \eqref{eq:B} and \eqref{eq:bsquare} follows that $B^p\equiv 1 \pmod{1+\alpha}, 1\le p\le \alpha$, so
\begin{equation}\label{eq:bigB}
1+\alpha\vert B^{\alpha}+B^{\alpha-1}+\cdots +B^2+B+1.
\end{equation}

Combining \eqref{eq:1} (for $\ell$) and \eqref{eq:bigB}, and taking into account \eqref{eq:account}, it follows that \eqref{eq:1} is true for $\ell+1$.
\end{proof}

\section{Proof of Proposition \ref{prop:1}} \label{sec:3-pre}

\begin{proof} The case $b=2$ is covered by  \cite[Theorem 12]{N}. If $b\ge 3$, choose $\alpha=b-2$ and apply Theorem \ref{thm:1}. We show now that the multipliers appearing in the proof of Theorem \ref{thm:1}, for a fixed base $b$, are all distinct. It follows from \eqref{eq:55} and \eqref{eq:66} that the multiplier for $N_k$ is given by:
\begin{equation}\label{eq:comp}
M=\frac{\frac{N_k}{2}}{s_b(N_k)}=\frac{\frac{b+\alpha}{2}\cdot\frac{b^{2(1+\alpha)^\ell}-1}{b^2-1}}{(1+\alpha)^{\ell+1}}.
\end{equation}

Note that $\alpha=b-2$. After algebraic manipulations, equation \eqref{eq:comp} becomes
\begin{equation*}
M=\frac{b^{2(1+\alpha)^{\ell}}-1}{(b-1)^{\ell}(b^2-1)}.
\end{equation*}

In order to show that the multipliers are distinct it is enough to show that the sequence of multipliers for $N_k$ is strictly increasing as a function of $\ell$, that is, we need to show that:
\begin{equation}\label{eq:22}
\frac{b^{2(1+\alpha)^{\ell}}-1}{(b-1)^{\ell}(b^2-1)}<\frac{b^{2(1+\alpha)^{\ell+1}}-1}{(b-1)^{\ell+1}(b^2-1)}.
\end{equation}

After algebraic manipulations \eqref{eq:22} becomes
\begin{equation}\label{eq:33}
(b-1)(b^{2(1+\alpha)^{\ell}}-1)< b^{2(1+\alpha)^{\ell+1}}-1.
\end{equation}
After denoting 
\begin{equation*}
B=b^{2(1+\alpha)^{\ell}}=b^{2(b-1)^{\ell}},
\end{equation*}
right hand side of \eqref{eq:33} factors as:
\begin{equation}\label{eq:131}
b^{2(1+\alpha)^{\ell+1}}-1=(b^{2(1+\alpha)^{\ell}}-1)(B^{\alpha}+B^{\alpha-1}+\cdots +B+1).
\end{equation}

Now \eqref{eq:33} follows from \eqref{eq:131} and the following inequality:
\begin{equation*}
b-1< b^{2(b-1)^{\ell}}, \ell\ge 0, \ell\ge 0, b\ge 3.
\end{equation*}
\end{proof}

\section{Proof of Theorem \ref{thm:new}}\label{sec:3}

\begin{proof} Let $b=p_1^{\alpha_1}p_2^{\alpha_2}\cdots p_k^{\alpha_k}, \alpha_i\ge 1$, $p_i$ prime, $1\le i\le k$. We recall that a base $b$ integer $N$ is divisible by $p_i^{\gamma}$ if the last $\gamma$  digits of $N$ form a base $b$ integer divisible by $p_i^{\gamma}$. 
Let $N=p_1^{\beta_1}p_2^{\beta_2}\cdots p_k^{\beta_k}w$, where $\gcd(w,b)=1$. Let $m=\max(\beta_1, \beta_2, \cdots , \beta_k)$. Let $L$ be the base $b$ integer equal to $p_1^{\beta_1}p_2^{\beta_2}\cdots p_k^{\beta_k}$. As $b\not\vert N$, the last digit of $L$ is not $0$.  Let $\ell$ be the length of $L$. Consider the base $b$ palindrome $P=[L^R(0)^{\land {m-\ell}}L]_b$, where $L^R$ is the reversal of base $b$-representation of $L$. As $P$ is divisible by $p_1^{\beta_1}p_2^{\beta_2}\cdots p_k^{\beta_k}$, this is the end of the proof if $w=1$. As $\gcd(w,b)=1$ Euler's theorem implies that $b^{\phi(w)}-1\pmod w$. 

Let $r$ be an integer divisible by $\phi(w)$ and greater than $l+m$, the length of $P$. Let $q\ge 1$ a multiple of $b^{\phi(w)}-1$. Consider
the infinite family of integers:
\begin{equation}\label{eq:newalpha}
\begin{gathered}
Q_{r,q}=\big[1\big((0)^{\land r-1}1\big)^{\land q}\big]_b=1+b^{r}+b^{2r}+\cdots+b^{qr}\\
=1+b^{r}+b^{2r}+\cdots+b^{qr}+q-q\\
=\big (b^{r}-1\big )+\big (b^{2r}-1\big )+\big (b^{3r}-1\big )+\cdots + \big (b^{qr}-1\big )+q.
\end{gathered}
\end{equation}
All terms in the last part of \eqref{eq:newalpha} are divisible by $b^{\phi(w)}-1$, so $Q_{r,q}$ is divisible by $b^{\phi(w)}-1$ and by $w$.
We finish the proof observing that $P\cdot Q_{r,q}$ is a base $b$ palindrome divisible by $N$.
\end{proof}

\section{Proof of Theorem \ref{thm:33}}\label{sec:44}

\begin{proof} Observe that:
\begin{equation}\label{eq:77}
\begin{gathered}
(b-1)\cdot (b-1)=b(b-2)+1=[(b-2)1]_b\\
(b-1)b^{k}+(b-1)b^{k}=b^k+(b-2)b^{k-1}=[1(b-2)0^{\land k}]_b.
\end{gathered}
\end{equation}

Using \eqref{eq:77} one has that

\begin{equation*}
\begin{gathered}
\ [1(b-1)]_b\cdot[(b-1)^{\land k}]_b=(b+b-1)\cdot \left ( \sum_{i=0}^{k-1}(b-1)b^i\right )\\
=\sum_{i=0}^{k-1} \Big ( (b-1)b^{i+1}+ \left ( b(b-2)+1\right )b^i\Big )\\
=\sum_{i=1}^{k}(b-1)b^i+\sum_{i=0}^{k-1} \left ( b(b-2)+1\right )b^i \\
=(b-1)b^{k}+\sum_{i=1}^{k-1}\Big ((b-1)+b(b-2)+1\Big ) b^i+b(b-2)+1\\
=(b-1)b^{k}+\sum_{i=1}^{k-1} (b-1)b^{i+1}+b(b-2)+1\\
=(b-1)b^{k}+(b-1)b^{k}+\sum_{i=1}^{k-2} (b-1)b^{i+1}+b(b-2)+1\\
=b^k+(b-2)b^{k-1}+\sum_{i=1}^{k-2} (b-1)b^{i+1}+b(b-2)+1\\
=[1(b-2)(b-1)^{\land k-2}(b-2)1]_b.
\end{gathered}
\end{equation*}
\end{proof} 

\section{Proof of Theorem \ref{thm:final}}\label{sec:final}

\begin{proof} a) Note that $s_b(N_k)=2a$. As $b$ is even, there exists an integer $M$ such that:
\begin{equation*}
2a\cdot M=[a(0)^{\land k+1}]_b.
\end{equation*}

The following computation shows that $N_k$ is a $b$-ARH number:
\begin{equation*}
\begin{gathered}
s_b(N_k)\cdot M+(s_b(N_k)\cdot M)^R\\
=[a(0)^{\land k+1}]_b+[a]_b=[a(0)^{\land k}a]_b=N_k.
\end{gathered}
\end{equation*}

To show that $N_k$ is not a $b$-Niven number observe that $N_k/a=[1(0)^{\land k}1]_b$ is not an odd number.

b) Note that $s_b(N_k)=2b$. As $b$ is even the multiplier $M=[(1(0)^{\land k})^{\land b}(0)^{\land kb}]_b/2$ is an integer. 

The following computation shows that $N_k$ is a $b$-ARH number:
\begin{equation*}
\begin{gathered}
s_b(N_k)\cdot M+(s_b(N_k)\cdot M)^R\\
=[(1(0)^{\land k})^{\land b}(0)^{\land kb+1}]_b+[((0)^{\land k}1)^{\land b}]_b=[(1(0)^{\land k})^{\land b}0 ((0)^{\land k}1)^{\land b}]_b=N_k.
\end{gathered}
\end{equation*}

To show that $N_k$ is not a $b$-Niven number observe that $N_k$ is not divisible by $b$.

c) The proof is similar to that of b).
\end{proof}

\section{Proof of Theorem \ref{thm:final2}}\label{sec:final2}

\begin{proof} a) Using that 
\begin{equation*}
(b-2)^2=b^2-4b+4=b(b-4)+4=[(b-4)4]_b,
\end{equation*}
an equivalent base b representation for $N_k$ is given by
\begin{equation*}
\begin{gathered}
N_k=[(b-2)(0)^{\land k-1}(b-4)5(0)^{\land k}(b-2)]_b, \text {if }b\not = 4,\\
N_k=[2(0)^{\land k-1}11(0)^{\land k}2]_4, \text {if }b = 4.
\end{gathered}
\end{equation*}
If $b\not = 4$ one has $s_b(N_k)=3(b-1)$ and if $b=4$ one has $s_4(N_k)=6$. To finish the proof of case a) it is enough to show that $\alpha_k$ is divisible by $s_b(N_k)$. 

If $b\not = 4$ one has that:
\begin{equation*}
\begin{gathered}
\alpha_k=b^{k+1}+b-2=b^{k+1}-1+b-1\\
=(b-1)\big (b^k+b^{k-1}+\dots +b^2+b+2 \big)
\end{gathered}
\end{equation*}
and
\begin{equation*}
\begin{gathered}
b^k+b^{k-1}+\dots +b^2+b+2\\
\equiv k+2 \pmod{3}\equiv 0\pmod{3},
\end{gathered}
\end{equation*}
where for the first congruence we use $b\equiv 1\pmod{3}$ and for the second congruence we use $k\equiv 1 \pmod{3}$.

If $b=4$, then clearly $\alpha_k$ is divisible by 2. Moreover
\begin{equation*}
\alpha_k=4^{k+1}+2=(3+1)^{k+1}+2\equiv 0\pmod{3},
\end{equation*}
which shows that $\alpha_k$ is divisible by 6.

b) Assume now $b=2$. Then an equivalent base 2 representation for $N_k$ is given by
\begin{equation*}
N_k=[1(0)^{\land k-1}10(0)^{\land k}1]_2,
\end{equation*}
so $s_2(N_k)=3$. To finish the proof we show that $\alpha_k$ is divisible by 3:
\begin{equation*}
\alpha_k=2^{k+1}+1=(3-1)^{k+1}+1\equiv 0\pmod{3}.
\end{equation*}
The congruence follows because $k$ is even.

\end{proof}

\section{Proof of Proposition \ref{prop:1-mult}} \label{sec:3-pre-mult}

\begin{proof} To prove the proposition we show that the multipliers from Theorem \ref{thm:final2} corresponding to various values of $k$ are distinct. This follows from the explicit formulas below, as all the sequences of multipliers are strictly increasing as functions of $k$:

If $b=2$ the sequence of multipliers is given by $M_k=\frac{2^{k+1}+1}{3}$.

If $b=4$  the sequence of multipliers is given by $M_k=\frac{4^{k+1}+2}{6}$.

If $b>4$ the sequence of multipliers is given by $M_k=\frac{b^{k+1}+b-2}{3(b-1)}$.

\end{proof}

\section{Proof of Theorem \ref{thm:last-thm}}\label{seq:last-thm}

\begin{proof} Let $I$ be a string of base $b$-digits. There exists an infinity of base $b$ strings $J$ such that $s_b([IJ]_b)$ is a power of $b$, say $b^k, k\ge 1$. Then the number $N_J=[IJ(0)^{\land k}]_b$ is a $b$-Niven number.
\end{proof}

\end{document}